\documentclass[11pt]{article}
\usepackage{mathrsfs}
\usepackage{stmaryrd}
\usepackage{amsfonts}
\usepackage{amsfonts,amsmath,amssymb,amscd}
\usepackage{shadow}
\usepackage{graphicx}
\usepackage{color}
\usepackage{longtable,colordvi}
\usepackage{longtable}
\allowdisplaybreaks
\newtheorem{theo}{Theorem}[section]
\newtheorem{defi}[theo]{Definition}
\newtheorem{lemm}[theo]{Lemma}
\newtheorem{prop}[theo]{Proposition}
\newtheorem{coro}[theo]{Corollary}
\newtheorem{remark}[theo]{Remark}

\def\qed{\hfill \rule{4pt}{7pt}}

\newcommand{\dotcup}{\ensuremath{\mathaccent\cdot\cup}}

\def\pf{\noindent {\it{Proof.} \hskip 2pt}}

\begin{document}
\begin{center}
{\LARGE\bf  The monomial basis and the $Q$-basis of the  Hopf
algebra of parking functions}
\end{center}

\begin{center}
Teresa X.S. Li

School of Mathematics and Statistics,\\
 Southwest University, Chongqing 400715, P.R. China

pmgb@swu.edu.cn

\end{center}

\begin{abstract}
Consider the vector space $\mathbb{K}\mathcal{P}$ spanned by  parking functions.
 By representing parking functions as labeled digraphs, Hivert, Novelli and Thibon
 constructed a cocommutative Hopf algebra PQSym$^{*}$ on $\mathbb{K}\mathcal{P}$.  The product and   coproduct of PQSym$^{*}$ are analogous to the product and coproduct of the Hopf algebra NCSym of symmetric
functions in noncommuting variables defined in terms of the power sum basis. In this paper, we view a parking function  as a word. We shall construct a Hopf algebra PFSym on $\mathbb{K}\mathcal{P}$  with a formal basis  $\{M_a\}$\ analogous to the monomial basis of  NCSym. By introducing a partial order on parking functions, we transform the basis  $\{M_a\}$ to another basis $\{Q_a\}$ via the M\"{o}bius inversion. We prove the freeness of PFSym by finding two free generating sets in terms of the $M$-basis and the $Q$-basis, and we show that PFSym is isomorphic to the Hopf algebra PQSym$^{*}$. It turns out that our construction, when restricted to permutations and non-increasing parking functions,  leads to a new way to approach the Grossman-Larson Hopf algebras of ordered trees and heap-ordered trees.
\end{abstract}


\section{Introduction}\label{section1}
The  algebra NCSym of symmetric functions in noncommuting variables
was first studied  by  Wolf \cite{Wolf1936} in 1936. Rosas and Sagan \cite{Rosas-Sagan2006} have found
 bases of NCSym analogous to the monomial, elementary, homogeneous and power sum bases of the algebra of symmetric functions in commuting variables. All these bases of NCSym are   indexed by set partitions. The Hopf algebra structure on NCSym was introduced by Bergeron, Reutenauer, Rosas and Zabrocki \cite{Bergeron-R-R-Z2008}.

Consider the vector space $\mathbb{K}\mathcal{P}$ spanned by parking functions. By representing parking functions as labeled digraphs, Hivert, Novelli and Thibon \cite{Hivert-Novelli-Thibon2008} introduced a Hopf algebra PQSym on $\mathbb{K}\mathcal{P}$ in their study of commutative  and cocommutative Hopf
algebras based on various combinatorial structures. The graded dual PQSym$^{*}$ of PQSym contains a Hopf subalgebra isomorphic to NCSym. This fact can be shown easily from  expressions for the product and the coproduct of NCSym in terms of the power sum basis \cite{Bergeron-H-R-Z2006}. To be more specific,  the product
of  PQSym$^{*}$ is given by shifted concatenation of parking functions and  the coproduct encodes all ways to divide connected components of a parking function into two parts and relabel each part. Similarly, in terms of the power sum basis, the product of NCSym is given by shifted union of set partitions and the coproduct is given by dividing all blocks of a partition into two parts and then standardize each part.
 In this paper, we are concerned with the problem of constructing a Hopf algebra on $\mathbb{K}\mathcal{P}$ with a basis  analogous to the monomial basis of  NCSym.  Using a formal basis $\{M_a\}$ indexed by parking functions, we  define a product $\star$ and a coproduct $\Delta$ on  $\mathbb{K}\mathcal{P}$  and we show that ${\rm PFSym}=(\mathbb{K}\mathcal{P},\star,\Delta)$  is a Hopf algebra  isomorphic to PQSym$^{*}$. As will be seen,  the basis $\{M_a\}$ of PFSym  is a natural generalization of the monomial basis of  NCSym.

To  define  the product  $\star$ and the coproduct $\Delta$ of PFSym, we view parking functions as words on positive integers. In this notation, a parking function
  can be decomposed into subwords via  its left-to-right minima. Then  the product $\star$ can be defined in terms of matchings between such decompositions and the coproduct $\Delta$ can be defined by dividing the decomposition into two  parts.  Also, based on this decomposition, we define the slash product and the split product  on parking functions, and then further introduce the concepts  of atomic parking functions and unsplitable parking functions. We shall show that the Hopf algebra PFSym is free by finding two free  generating sets. The first one consists of those basis elements $M_{\alpha}$ indexed  by  unsplitable parking functions. The second one, indexed by atomic parking functions, consists of elements  from another basis $\{Q_a\}$.   The $Q$-basis  for PFSym is related to the  basis $\{M_a\}$ via the M\"{o}bius inversion on a partial order on parking functions.  This basis is a natural analog of the $q$-basis of the algebra NCSym introduced by Bergeron and Zabrocki \cite{Bergeron-Zabrocki2005}.

This paper is organized as follows. In Section \ref{sec:2}, we give an overview
 of the Hopf  algebra ${\rm NCSym}$ of symmetric functions in noncommuting variables. The product of ${\rm NCSym}$  will be described in terms of matchings between set partitions.  In Section \ref{sec:3}, we introduce the   $LR$-decomposition of a parking function. By means of this decomposition and a formal basis $\{M_a\}$ indexed by parking functions, we define a product $\star$ and a coproduct $\Delta$ on $\mathbb{K}\mathcal{P}$  and we show that ${\rm PFSym}=(\mathbb{K}\mathcal{P},\star,\Delta)$ is a Hopf algebra.  In Section \ref{sec:4}, we introduce another basis $\{Q_a\}$ via a partial order on parking functions and show  the freeness of ${\rm PFSym}$ in terms of the $M$-basis and the $Q$-basis.  In the last section, we consider Hopf subalgebras of ${\rm PFSym}$. In addition to a Hopf subalgebra isomorphic to NCSym, we find  Hopf subalgebras isomorphic to the Grossman-Larson Hopf algebras of  ordered trees and heap-ordered trees \cite{Grossman-Larson1989}.

\section{The Hopf algebra ${\rm NCSym}$}\label{sec:2}

In this section, we  give an overview of the Hopf algebra ${\rm NCSym}$ of symmetric functions in noncommuting variables. This Hopf algebra was introduced by   Bergeron, Reutenauer, Rosas and Zabrocki \cite{Bergeron-R-R-Z2008} in their study of the connections between the algebra of symmetric functions in commuting variables and the algebra of symmetric functions in noncommuting variables.
Instead of employing the lattice structure of partitions, we
use matchings between partitions to  describe the product of ${\rm NCSym}$
 with respect to the monomial basis. As will be seen, this description is of importance in providing us ideas to construct the Hopf algebra PFSym of parking functions.

 Throughout this paper,  $\mathbb{K}$ stands for a field of characteristic zero and  $\mathbb{K}S$ denotes the vector space with the set $S$ as basis. For the background on Hopf algebras, see  Abe \cite{Abe1980}  and Sweedler \cite{Sweedler1969}.

The Hopf algebra ${\rm NCSym}$ can be formally described in terms of
partitions. Let $[n]=\{1,2,\ldots,n\}.$  A partition of $[n]$ is a set $\{B_1,B_2,\ldots,B_k\}$
of pairwise disjoint nonempty subsets of $[n]$ whose union is $[n]$. The subsets $B_i$
are called blocks of $\pi$. Without loss of generality, we may assume that the blocks of
a partition are arranged in the increasing order of their minimal elements and  the
elements in each block are written in increasing order. Denote by $\Pi_n$ the set of
 partitions of $[n]$. Set $\Pi_0=\{\emptyset\}$ and $\Pi=\cup_{n\geq 0}\Pi_n$. Let
$$
{\rm NCSym}=\bigoplus_{n\geq 0}\mathbb{K}\{M_\pi:\pi\in \Pi_n\},
$$
where $\{M_\pi\}$ is a formal basis called monomial basis and
the space $\mathbb{K}\{M_\emptyset\}$ is regarded as $\mathbb{K}$ by identifying
$M_{\emptyset}$ with $1$.

The  product of ${\rm NCSym}$ is defined as follows.  Firstly, for any $\pi\in\Pi,$ we set $M_{\pi} M_{\emptyset}=M_{\emptyset}M_{\pi}=M_{\pi}$.
Suppose that $m,n\geq 1$, $\pi=\{B_1,B_2,\ldots,B_r\}\in\Pi_m$ and  $\sigma=\{C_1,C_2,\ldots,C_s\}\in\Pi_n$. Let $\sigma+m=\{C_1+m,C_2+m,\ldots,C_s+m\},$ where $C_i+m$ denotes the set obtained by adding $m$ to each element in $C_i$.
We shall consider the set $R(\pi,\sigma)$ of all  matchings between $\pi$ and $\sigma+m$, where the blocks of $\pi$ and $\sigma+m$ are considered as
 vertices. Here a matching between two vertex sets $X$ and $Y$ means a bipartite graph with bipartition $(X, Y)$ such that there are no two edges
 having a common vertex. For example, if $\pi=\{B_1, B_2, B_3\}$ and $\sigma=\{C_1, C_2, C_3, C_4\}$, the following diagram represents a matching with
 two edges $(B_1, C_3+m)$ and $(B_2, C_1+m)$.
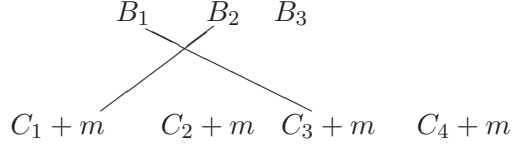
 \begin{center}\setlength{\unitlength}{1mm}
\begin{picture}(20,18)
\put(-12,15){$B_1$}\put(0,15){$B_2$}\put(9,15){$B_3$}
\put(-26,0){$C_1+m$}\put(-6,0){$C_2+m$}\put(10,0){$C_3+m$}
\put(28,0){$C_4+m$}
\put(-14,3){\line(4,3){15}}
\put(14,3){\line(-2,1){22}}
\end{picture}
\end{center}

  Given a matching in $R(\pi, \sigma)$,
  we can construct a partition  from the blocks of
 $\pi$ and $\sigma+m$ by  combining  blocks $B_i$ and $C_j+m$  if they form an edge of the matching. Let $S(\pi, \sigma)$ denote
 the set of such partitions obtained from matchings in $R(\pi, \sigma)$.
Then the product  of ${\rm NCSym}$ is defined  by
\begin{align}\label{product of NCSym}
M_\pi M_{\sigma}=\sum_{\tau\in S(\pi,\sigma)}M_{\tau}.
\end{align}
For example, we have
\begin{align*}
M_{\{\{1,3\},\{2,4\}\}} M_{\{\{1,3,4\},\{2\}\}}=&M_{\{\{1,3\},\{2,4\},\{5,7,8\},\{6\}\}}+M_{\{\{1,3,5,7,8\},\{2,4\},\{6\}\}}\\[5pt]
&+M_{\{\{1,3\},\{2,4,5,7,8\},\{6\}\}}
+M_{\{\{1,3,6\},\{2,4\},\{5,7,8\}\}}\\[5pt]
&+M_{\{\{1,3\},\{2,4,6\},\{5,7,8\}\}}
+M_{\{\{1,3,5,7,8\},\{2,4,6\}\}}\\[5pt]
&+M_{\{\{1,3,6\},\{2,4,5,7,8\}\}}.
\end{align*}

To define the coproduct on ${\rm NCSym}$, we need the notion of standardization.
Suppose that $\tau=\{B_1,B_2,\ldots,B_k\}$ is a family of disjoint nonempty finite sets of integers.  Define the standardization of $\tau$, denoted
by ${\rm st}(\tau)$, to be the set partition obtained from $\tau$ by substituting
the smallest element by $1$, the second smallest element by $2$, and so on. By convention, we let ${\rm st}(\emptyset)=\emptyset$. For example, let $\tau=\{\{2,5,7\},\{4,8\},\{6,9\}\}$. Then we have  ${\rm st}(\tau)=\{\{1,3,5\},\{2,6\},\{4,7\}\}$.

The coproduct
\[
\Delta\colon  {\rm NCSym}\longrightarrow {\rm NCSym}\otimes {\rm NCSym}
\]
is defined by
\begin{align}\label{coproduct of NCSym}
\Delta(M_{\pi})=\sum_{\pi_1\dotcup\pi_2=\pi}M_{{\rm st}(\pi_1)}\otimes M_{{\rm st}(\pi_2)},
\end{align}
where by the notation
 \[
\sum_{\pi_1\dotcup\pi_2=\pi}
 \] we mean that the sum ranges over  ordered pairs $(\pi_1,\pi_2)$ of disjoint subsets of $\pi$ such  that $\pi_1\cup\pi_2=\pi$, see \cite{Lauve-Mastnak2011}.

For example, we have
\begin{align*}
\Delta(M_{\{\{1,4,6\},\{2,5\},\{3\}\}})
=&1\otimes M_{\{\{1,4,6\},\{2,5\},\{3\}\}}
+M_{\{\{1,2,3\}\}}\otimes M_{\{\{1,3\},\{2\}\}}\\[5pt]
&+M_{\{\{1,2\}\}}\otimes M_{\{\{1,3,4\},\{2\}\}}+M_{\{\{1\}\}}\otimes M_{\{\{1,3,5\},\{2,4\}\}}\\[5pt]
&+M_{\{\{1,3\},\{2\}\}}\otimes M_{\{\{1,2,3\}\}}+M_{\{\{1,3,4\},\{2\}\}}\otimes M_{\{\{1,2\}\}}\\[5pt]
&+M_{\{\{1,3,5\},\{2,4\}\}}\otimes M_{\{\{1\}\}}
+M_{\{\{1,4,6\},\{2,5\},\{3\}\}}\otimes 1.
\end{align*}

Equipped with the product  and coproduct defined as above,
${\rm NCSym}$ is  a graded cocommutative Hopf
algebra with  unit $M_{\emptyset}=1\in\mathbb{K}$ and counit given by
\begin{align*}
\varepsilon(M_{\pi})=\left\{
                       \begin{array}{ll}
                         1, & \hbox{if $\pi=\emptyset$;} \\[5pt]
                         0, & \hbox{otherwise.}
                       \end{array}
                     \right.
\end{align*}

\section{The Hopf algebra PFSym of parking functions}\label{sec:3}

In this section, based on a formal basis analogous to the monomial basis of ${\rm NCSym}$, we define a product $\star$ and a coproduct $\Delta$ on the vector space $\mathbb{K}\mathcal{P}$ spanned by parking functions. We show that, with these operations, $\mathbb{K}\mathcal{P}$  is indeed a Hopf algebra.

We begin by introducing some basic notation and terminology. Let $[n]^k$ denote the set of words of length $k$ on $[n]$.
A word $a=a_1a_2\cdots a_n\in [n]^n$ is called a parking function
if its nondecreasing rearrangement $b_1b_2\cdots b_n$ satisfies $b_i\leq i$
for $1\leq i\leq n$. Let $\epsilon$ denote the unique empty word of length $0$. Denote by $P_n$ the set
of  parking functions of length $n$ and set $P=\cup_{n\geq 0}P_n$, where
 $P_0=\{\epsilon\}$.

Suppose that  $a=a_1a_2\cdots a_n$ is a word on positive integers and $1\leq i\leq n$. We say that $i$ is a position of left-to-right minimum of $a$ if $i=1$ or $i\geq 2$ and  $a_i<a_j$ for any $j< i$. Let $lr(a)=\{i_1,i_2,\ldots,i_s\}$ denote the  set of positions of left-to-right minima of $a$, where  $i_1<i_2<\cdots< i_s$. Set $i_{s+1}=n+1$.
For $1\leq j\leq s$, let
   \[
   w_j=a_{i_{s+1-j}}a_{i_{s+1-j}+1}\cdots a_{i_{s+2-j}-1}.
   \]
We call the sequence $(w_1,w_2,\ldots,w_s)$ of words  the LR-decomposition of $a$, denoted by $F_a$.

For example, let $a=56357622315\in P_{11}$. Then we have $lr(a)=\{1,3,7,10\}$ and
 \[F_a=(15,223,3576,56).\]

\begin{remark} \label{remark for F_a}
 For a word $a$, let ${\rm min}(a)$ denote the minimal integer appearing in  $a$. A word $a=a_1a_2\ldots a_k$ is said to be dominated if $a_1={\rm min}(a)$.  Let
 \[
 \mathcal{F}_n=\{F_a\ |\ a\in P_n\},\ n\geq 1
 \] and set $\mathcal{F}_0=\{(\,)\}=\{F_{\epsilon}\}$.
It is easy to verify that a sequence $F=(w_1,w_2,\ldots,w_k)\in\mathcal{F}_n$
if and only if
\begin{itemize}
\item[(1)] the  word $w_k\cdots w_2\cdot w_1$ obtained by concatenating $w_i$ is a parking function;
\item[(2)] for  $1\leq i\leq k$, the word  $w_i$ is dominated;
\item[(3)] ${\rm min}(w_1)<{\rm min}(w_2)<\cdots<{\rm min}(w_k)$.
\end{itemize}
Moreover, for each $n \geq 0$, $a\mapsto F_a$ is a bijection between $P_n$ and $\mathcal{F}_n$.
\end{remark}

In the following, unless otherwise specified, we always identify  a sequence
$F=(w_1,w_2,\ldots,w_k)$ of dominated  words   such that
${\rm min}(w_1)<{\rm min}(w_2)<\cdots<{\rm min}(w_k)$ with the underlying set
$\{w_1,w_2,\ldots,w_k\}$, since we can assume that the words in the set are arranged
in increasing order of their minimal elements.   So  the notations such as
$w_i\in F$, $F\cap G$, $F\cup G$, and $F\setminus G$ make sense.

Let
$$
\mathbb{K}\mathcal{P}=\bigoplus_{n\geq 0}\mathbb{K}\{M_a:a\in P_n\},
$$ where $M_a$ denotes a formal basis element indexed by a parking function $a$.  We call the basis $\{M_a\ |\ a\in P\}$  monomial basis of $\mathbb{K}\mathcal{P}$ and  regard  $\mathbb{K}\{M_\epsilon\}$ as $\mathbb{K}$ by identifying $M_{\epsilon}$ with $1$.

The product on $\mathbb{K}\mathcal{P}$ is a natural generalization of that in ${\rm NCSym}$, as described in Section 2. Firstly, set $M_{a}\star M_{\epsilon}=M_{\epsilon}\star M_{a}=M_{a}$ for any parking function $a$. Now suppose that $m\geq 1, n\geq 1$, $a=a_1a_2\cdots a_m\in P_m$ and
$b=b_1b_2\cdots b_n\in P_n$. Let
\[F_a=(u_1,u_2,\ldots,u_r),\ \
F_b=(w_1,w_2,\ldots,w_s)\]
and $$F_b+m=(w_1+m,w_2+m,\ldots,w_s+m),$$
where $w_i+m$ denotes the word obtained by adding $m$ to each integer in $w_i$. Let $R(a,b)$ denote the set of matchings between $F_a$ and
 $F_b+m$. Given a matching $\Theta$ in $R(a,b)$, we  define a set  $F_a\Theta F_b$ of words from
 $F_a$ and $F_b+m$ by  concatenating words $u_i$ and $w_j+m$  if they form an edge in $\Theta$. Obviously, we have $F_a\Theta F_b\in\mathcal{F}_{m+n}$ and so we get a unique parking function, denoted by $a\Theta b$, whose $LR-$decomposition is $F_a\Theta F_b$.
Define the product  of $M_a$ and $M_b$  by
\begin{align}\label{product of KP}
M_a\star M_b=\sum_{\Theta\in R(a,b)}M_{a\Theta b}.
\end{align}
For example, let \[ a=211, \ b=353112.\]
Then we have $$F_a=(11,2),\ F_b+3=(445,686)$$ and
\begin{align*}
\{F_a\Theta F_b\colon \Theta\in R(a,b)\}
&=\{(11,2,445,686),(11,2686,445),(11686,2,445),\\[5pt]
&(11,2445,686),(11445,2,686),(11445,2686),(11686,2445)\}.
\end{align*}
So we have
\begin{align*}
M_a\star M_b=&M_{686445211}+M_{445268611}+M_{445211686}+M_{686244511}\\[5pt]
&+M_{686211445}+M_{268611445}+M_{244511686}.
\end{align*}

Now we proceed to the definition of  the coproduct map
\[
\Delta\colon \mathbb{K}\mathcal{P}\longrightarrow \mathbb{K}\mathcal{P}\otimes\mathbb{K}\mathcal{P}.
\]
To this end, we shall use the parkization ${\rm Park}(a)$ of a word $a$, which was used by J.-C. Novelli and J.-Y. Thibon in \cite{Novelli-Thibon2007}. For completeness, let us give an overview of the definition of ${\rm Park}(a)$. Suppose that $a=a_1a_2\cdots a_n$
is a word on positive integers. Set
\[
d(a)= {\rm min}\{i: \ |\{j:a_j\leq i\}| < i\}.
\]
If $d(a)=n+1$, then we find that $a$ is a parking function and let ${\rm Park}(a)=a$.
 Otherwise, let $\bar{a}$ be the word obtained from $a$ by decrementing  all the elements greater than  $d(a)$ and let ${\rm Park}(a)={\rm Park}(\bar{a})$.   Since $\bar{a}$ is smaller than $a$ in the
lexicographic
order, the algorithm terminates and always returns a parking function ${\rm Park}(a),$
called the parkization of $a$. It can be seen that when $a$ is a word without repeated letters,
${\rm Park}(a)$ coincides with the standardization of $a$, which is obtained from $a$ by substituting the smallest element by $1$, the second smallest element by $2$, and so on.  Set ${\rm Park}(\epsilon)=\epsilon$.

For example, we give the algorithm of computing ${\rm Park}(875221)$ as follows:
\begin{align*}
a&=875221 \qquad  d(a)=4\\[3pt]
\bar{a}&=764221 \qquad  d(\bar{a})=5\\[3pt]
\bar{\bar{a}}&=654221 \qquad  d(\bar{\bar{a}})=7,
\end{align*}
so we have ${\rm Park}(875221)={\rm Park}(764221)={\rm Park}(654221)=654221$.


Now the coproduct $\Delta$ is defined by
\begin{align}\label{coproduct of KP}
\Delta(M_{a})=\sum_{F_{a'}\dotcup F_{a''}=F_a}M_{{\rm Park}(a')}\otimes M_{{\rm Park}(a'')},
\end{align}
where the sum ranges over  ordered pairs $(a',a'')$ of subwords of $a$ such  that $F_{a'}\cap F_{a''}=\emptyset$ and $F_{a'}\cup F_{a''}=F_a$.

 For example, let $a=445132$, then we have
$F_a=(132,445)$ and
\begin{align*}
\Delta(M_a)=&M_{\epsilon}\otimes M_{445132}+M_{{\rm Park}(445)}\otimes M_{{\rm Park}(132)}\\[5pt]
&+M_{{\rm Park}(132)}\otimes M_{{\rm Park}(445)}+M_{445132}\otimes M_{\epsilon}\\[5pt]
=&1\otimes M_{445132}+M_{112}\otimes M_{132}+M_{132}\otimes M_{112}+M_{445132}\otimes 1.
\end{align*}

Let ${\rm PFSym}=(\mathbb{K}\mathcal{P},\star,\Delta)$ denote the vector space $\mathbb{K}\mathcal{P}$ equipped with the product $\star$ and the coproduct $\Delta$.

\begin{theo}\label{Main theorem}
${\rm PFSym}$ is a cocommutative Hopf algebra with unit $M_{\epsilon}$ and counit given by
\begin{align}\label{counit of KP}
\varepsilon(M_{a})=\left\{
                     \begin{array}{ll}
                       1, & \hbox{if $a=\epsilon$;} \\[5pt]
                       0, & \hbox{otherwise.}
                     \end{array}
                   \right.
\end{align}
\end{theo}

We will prove Theorem \ref{Main theorem} through the following  steps.

\begin{prop}\label{algebra}
With the product $\star$ defined by \eqref{product of KP}, $\mathbb{K}\mathcal{P}$
is an associative algebra.
\end{prop}
\pf Suppose that  $a\in P_l$, $b\in P_m$ and $c\in P_n$. Let
$F_a=(u_1,u_2,\ldots,u_r)$, $F_b=(v_1,v_2,\ldots, v_s)$ and $F_c=(w_1,w_2,\ldots, w_t)$.
We need to show that
\begin{align}\label{associativity}
(M_a\star M_b)\star M_c=M_a\star (M_b\star M_c).
\end{align}
 By definition, we have
\begin{align*}
(M_a\star M_b)\star M_c=\sum_{\Theta\in R(a,b)}\ \sum_{\Lambda\in R(a\Theta b,c)}
M_{(a\Theta b)\Lambda c}.
\end{align*}
and
\begin{align*}
M_a\star (M_b\star M_c)=\sum_{\Theta'\in R(b,c)}\
\sum_{\Lambda'\in R(a,b\Theta' c)}M_{a\Lambda'(b\Theta' c)}.
\end{align*}
To prove \eqref{associativity}, it suffices to find a bijection
\begin{align*}
 \{(\Theta,\Lambda)\ |\ \Theta\in R(a,b), \Lambda\in R(a\Theta b,c)\}&\longrightarrow
\{(\Theta',\Lambda')\ |\ \Theta'\in R(b,c), \Lambda'\in R(a, b\Theta' c) \}\\[5pt]
(\Theta,\Lambda)&\longmapsto (\Theta',\Lambda')
\end{align*}
 such that $(a\Theta b)\Lambda c=a\Lambda'(b\Theta' c).$ Suppose that $\Theta\in R(a,b)$ and $\Lambda\in R(a\Theta b,c)$. Let $\Theta'$ be the matching between $F_b$ and $F_c+m$ whose edges are those $(v_j,w_k+m)$ such that $v_j\in F_b$, $w_k\in F_c$ and
 \begin{itemize}
 \item[(i)] $(v_j+l,w_k+l+m)$ is an edge in $\Lambda$ or
 \item[(ii)] $(u_i\cdot (v_j+l),w_k+l+m)$ is an edge in $\Lambda$ for some $u_i\in F_a$.
 \end{itemize}
And we let $\Lambda'$ be the matching between $F_a$ and $F_{b\Theta' c}+l$ that consists of the following three kinds of edges:
\begin{itemize}
\item[(i)] $(u_i,v_j+l)$, where  $u_i\in F_a$, $v_j\in F_b$ and $(u_i,v_j+l)$ is an edge in $\Theta$;
\item[(ii)] $(u_i, (v_j+l)\cdot (w_k+l+m))$, where  $u_i\in F_a$, $v_j\in F_b$, $w_k\in F_c$ and $(u_i\cdot (v_j+l),w_k+l+m)$ is an edge in $\Lambda$;
\item[(iii)] $(u_i,w_k+l+m)$, where $u_i\in F_a$, $w_k\in F_c$ and $(u_i,w_{k}+l+m)$ is an edge in $\Lambda$.
\end{itemize}
Then it can be checked easily that $(\Theta,\Lambda)\mapsto (\Theta',\Lambda')$
 is the desired  bijection.
  The proof is complete. \qed

By the definition
of the coproduct  and the counit given in \eqref{coproduct of KP} and
\eqref{counit of KP}, we can  verify the
coassociativity and the counitary property without any difficult, so we obtain

\begin{prop}\label{coalgebra}
With the coproduct $\Delta$ defined by \eqref{coproduct of KP} and counit $\varepsilon$ defined by \eqref{counit of KP},
$\mathbb{K}\mathcal{P}$ is a
cocommutative coalgebra.
\end{prop}

By \eqref{counit of KP}, it is easy to see that the counit is an algebra map.
So it suffices to show that the coproduct map is also an  algebra  map. To this end, we need two lemmas. We begin with mentioning a simple observation and introducing some notation.

Firstly, a simple observation shows that, for any word $a$, the relative order
of integers in $a$ is invariant under the operator  ${\rm Park}$. So there is a one-to-one correspondence $\iota: F_{{\rm Park}(a)}\rightarrow F_a$ such that the relative order of integers in each word of $\iota(u)$ is the same as that of $u$. In fact, if $\bar{a}$ is the word appearing in the definition of ${\rm Park}(a)$, then there is a similar bijection between $F_{\bar{a}}$ and $F_{a}$. We will use this fact occasionally in the following.

Suppose that $a, b$ are two words on positive integers. For an integer $k$ and a matching $\Theta_k$ between $F_a$ and $F_b+k$, let $a\Theta_k b$ denote the word whose LR-decomposition is obtained from $F_a$ and $F_b+k$ by concatenating $u$ and $w+k$ if $(u,w+k)$ is an edge in $\Theta$.  Note that, when $a$ and $b$ are parking functions, we use the simpler notation $a\Theta b$ instead of $a\Theta_{l(a)} b$, where $l(a)$ denotes the length of $a$.

Now assume that $a\in P_m$, $b\in P_n$, $\Theta\in R(a,b)$. Suppose that $a'$ and $b'$ are subwords of $a$ and $b$ respectively such that
\[F_{a'}\subseteq F_{a},\ \ F_{b'}\subseteq F_{b},\ \ l(a')=m'\ \ {\rm and}\ \ l(b')=n'.
\]
 By restricting  $\Theta$ to $F_{a'}$ and $F_{b'}+m$, we get a matching $\Theta_m$ between $F_{a'}$ and $F_{b'}+m$, so the notation $a'\Theta_m b'$ makes sense. Let $\Theta'$ be the matching between $F_{{\rm Park}(a')}$ and $F_{{\rm Park}(b')}+m'$ with edge set
\begin{align*}
\{(u,w+m') \ |\  \iota_1(u)\cdot(\iota_2(w)+m)\in F_{a'\Theta_{m}b'}\},
\end{align*}
where $\iota_1: F_{{\rm Park}(a')}\rightarrow F_{a'}$ and $\iota_2: F_{{\rm Park}(b')}\rightarrow F_{b'} $ denote the natural bijections mentioned above. Note that $\Theta'$ only depends on $\Theta_m$.


\begin{lemm}\label{Parkized lemm}We have
$${\rm Park}(a'\Theta_{m}b')={\rm Park}(a')\Theta'{\rm Park}(b').$$
\end{lemm}
\pf Firstly, we show the following two claims.

Claim (1):  Suppose that $c\in P_{m'}$ and $b=b_1b_2\cdots b_k$ is any word on positive integers. If $m\geq m'$, then for any matching $\Theta_m$ between $F_c$ and $F_{b}+m$, we have
 \[
{\rm Park}(c\,\Theta_m b)=c\,\Theta' {\rm Park}(b)
\]
where $\Theta'$ is the matching between $F_c$ and $F_{{\rm Park}(b)}+m'$ with edge set
\[
\{(u,w+m')\ |\ u\cdot(\iota(w)+m)\in F_{c\Theta_m b}\}
\]
and $\iota$ is the natural bijection between $F_{{\rm Park}(b)}$ and $F_{b}$.
 We proceed by double induction on $m$ and $\sum_{j=1}^{k}b_{j}$.  It is quite easy to verify the claim for $m=m'$ and $b'=11\cdots1$.

Since $c\in P_{m'}$, we have
$d(c\,\Theta_m b)\geq m'+1$, where $d(x)$ is the integer appearing in the definition of ${\rm Park}(x)$. If $d(c\,\Theta_m b)>m'+1$, then  there exists some $b_j$ such that $b_{j}+m\leq m'+1$. It follows that $m=m'$ and $b_{j}=1.$
In this case, it is easy to see that $d(c\Theta_m b)=d(b)+m$ and
\begin{align*}
{\rm Park}(c\,\Theta_m b)=\left\{
                            \begin{array}{ll}
                              {\rm Park}(c\, \bar{\Theta}_m \bar{b} ), & \hbox{if $b\notin P_{k}$;} \\
                              c\, \Theta_m b, & \hbox{if $b\in P_{k}$}
                            \end{array}
                          \right.
\end{align*}
where $\bar{\Theta}_m$ is the matching between $F_{c}$ and $F_{\bar{b}}+m$ with edge set
\[
\{(u,w+m)\ |\ u\cdot(\bar{\iota}(w)+m)\in F_{c\Theta_m b}\}
\]
and $\bar{\iota}$ is the natural bijection between $F_{\bar{b}}$ and $F_{b}$. If $b\in P_{k}$, we are done. Otherwise, we can complete the proof by induction since $\sum_{j=1}^{k}\bar{b}_{j}<\sum_{j=1}^{k}b_{j}$ and $\bar{\Theta}'=\Theta'$, where $\bar{\Theta}'$ is the matching between $F_c$ and $F_{{\rm Park}(\bar{b})}+m'$ obtained from $\bar{\Theta}_{m}$ in the same way as that $\Theta'$ is obtained from $\Theta_m$.

If $d(c\Theta_m b)=m'+1$, then $m'<m$ or $m=m'$ and $b_{i}>1$ for all $1\leq i\leq k$. If $m'<m$, by the definition of parkization, we have
\[
{\rm Park}(c\,\Theta_m b)={\rm Park}(c\,\bar{\Theta}_{m-1} b),
\]
where $\bar{\Theta}_{m-1}$ denotes the matching between $F_{c}$ and $F_{b}+(m-1)$ with edge set
 \[
\{(u,w+(m-1))\ |\ u\cdot (w+m)\in F_{c\,\Theta_m b}\}.
\]Therefore, we can complete the proof by induction since $m'\leq m-1$ and $\bar{\Theta}'=\Theta'$. If $m=m'$ and $b_{i}>1$ for all $1\leq i\leq k$, then by the definition of parkization we have
\[
{\rm Park}(c\,\Theta_m b)={\rm Park}(c\,\bar{\Theta}_m (b-1)),
\]
where $b-1$ denotes the word obtained from $b$ by subtracting $1$ from each element and $\bar{\Theta}_{m}$ denotes the matching between $F_{c}$ and $F_{b-1}+m$ with edge set
\[
\{(u,w+m)\ |\ u\cdot ((w+1)+m)\in F_{c\Theta_m b}\}.
\]
Again, we have $\bar{\Theta}'=\Theta'$.  Since $\sum_{j=1}^{k}(b_{j}-1)<\sum_{j=1}^{k}b_j$,  by induction, we are done.

Claim (2):  Suppose that $c\in P_m$ and $c'=c_{i_1}c_{i_2}\cdots c_{i_{m'}}$ is any subword of $c$. Let $j(c')$ be the number of steps involved in the algorithm of computing ${\rm Park}(c')$ , namely,
    \begin{align*}
    j(c')=\left\{
            \begin{array}{ll}
              0, & \hbox{if $c'\in P_{m'}$;} \\
              j(\bar{c'})+1, & \hbox{otherwise,}
            \end{array}
          \right.
    \end{align*}
where $\bar{c'}$ denotes the word appearing in the definition of ${\rm Park}(c')$. Then we have $j(c')\leq m-m'.$

The claim is trivially true if $c'\in P_{m'}$. Assume that $c'\notin P_{m'}$. Then $m'<m$ and $j(c')=j(\bar{c'})+1$ . Let $c^{*}$ be the word obtained from $c$ by replacing each $c_{i_k}$ by $\bar{c'}_k$ and then deleting a maximal element in the complementary subword $c_{j_1}c_{j_2}\cdots c_{j_{m-m'}},$ where $\{j_1<j_2<\cdots<j_{m-m'}\}=[m]\setminus \{i_1,i_2,\ldots,i_{m'}\}$. Clearly, $c^{*}\in P_{m-1}$ and $\bar{c'}$ is a subword of $c^{*}$. By an inductive argument on $m$, we get that
\[
j(c')=j(\bar{c}')+1\leq m-1-m'+1=m-m'.
\]

Now we are ready to prove the lemma. If $a'\in P_{m'}$, by Claim (1), we are done. Assume that $a'\notin P_{m'}$, then $m'<m$ and $d(a'\Theta_m b')\leq m'$. By the definition of parkization, we have
\[
{\rm Park}(a'\Theta_m b')={\rm Park}({\rm Park}(a')\,\bar{\Theta}_{m-j(a')}\, b')
\]
where $j(a')$ denotes the number of steps involved in the algorithm of computing ${\rm Park}(a')$ and  $\bar{\Theta}_{m-j(a')}$ denotes the matching between $F_{{\rm Park}(a')}$ and $F_{b'}+(m-j(a'))$ with edge set
\begin{align*}
\{(u,v+(m-j(a')))\ |\ \iota_1(u)\cdot(v+m)\in F_{a'\Theta_m b'}\},\ \  \ \ (\iota_1: F_{{\rm Park}(a')} \rightarrow F_{a'}).
\end{align*}
Clearly, by the construction of $\Theta'$ , we have  $\bar{\Theta}'=\Theta'$. By Claim (2), we have $m'\leq m-j(a')$.  Then the proof follows from  Claim (1).
\qed

\begin{lemm}\label{main lemm}
Suppose that $a\in P_m$ and $b\in P_n$. Let
\[
X(a,b)=\{(\Theta,c',c'')\ |\ \Theta\in R(a,b), F_{c'}\dotcup F_{c''}= F_{a\Theta b}\}
\] and
\begin{align*}
Y(a,b)=\{&(a',a'',b',b'', \Theta', \Theta'')\ |\ F_{a'}\dotcup F_{a''}=F_a, F_{b'}\dotcup F_{b''}=F_b, \\[5pt]
&\Theta'\in R({\rm Park}(a'),{\rm Park}(b')), \Theta''\in R({\rm Park}(a''),{\rm Park}(b''))\}.
\end{align*}
Then there exists a bijection
\begin{align*}
X(a,b)&\longrightarrow Y(a,b)\\[5pt]
(\Theta,c',c'')&\longmapsto (a',a'',b',b'', \Theta', \Theta'')
\end{align*}
such that
\begin{align}\label{fact1}
{\rm Park}(c')={\rm Park}(a')\Theta'{\rm Park}(b')
\end{align}
 and
\begin{align}\label{fact2}
{\rm Park}(c'')={\rm Park}(a'')\Theta''{\rm Park}(b'').
\end{align}
\end{lemm}
\pf Suppose that $F_a=(u_1,u_2\ldots,u_r)$,
$F_b=(w_1,w_2,\ldots, w_s)$, and $(\Theta,c',c'')\in X(a,b)$. Let $a'$, $a''$, $b'$ and $b''$ be the words  such that
\begin{align*}
&F_{a'}=\{u\in F_a\ |\ u\in F_{c'}\ {\rm or}\ u\cdot(w+m)\in F_{c'}\ {\rm for\ some\ }w\in F_b\}\\[5pt]
&F_{a''}=\{u\in F_a\ |\ u\in F_{c''}\ {\rm or}\ u\cdot(w+m)\in F_{c''}\ {\rm for\ some\ }w\in F_b\}\\[5pt]
&F_{b'}=\{w\in F_b\ |\ w+m\in F_{c'}\ {\rm or}\ u\cdot(w+m)\in F_{c'}\ {\rm for\ some}\ u\in F_a\}\\[5pt]
&F_{b''}=\{w\in F_b\ |\ w+m\in F_{c''}\ {\rm or}\ u\cdot(w+m)\in F_{c''}\ {\rm for\ some}\ u\in F_a\}.
\end{align*}
Then clearly we have $F_{a'}\dotcup F_{a''}=F_a$ and $F_{b'}\dotcup F_{b''}=F_b$. Assume that ${\rm Park}(a')\in P_{m'}$, then ${\rm Park(a'')}\in P_{m-m'}$. Now let
$\Theta'$ be the matching between $F_{{\rm Park}(a')}$ and $F_{{\rm Park}(b')}+m'$
with edge set:
\[
\{(u',w'+m')\ |\ \iota_{1}(u')\cdot(\iota_{2}(w')+m)\in F_{c'}\},
\]
and let $\Theta''$ be the matching between $F_{{\rm Park}(a'')}$ and $F_{{\rm Park}(b'')}+m-m'$
with edge set:
\[
\{(u'',w''+m-m')\ |\ \iota_{3}(u'')\cdot(\iota_{4}(w'')+m)\in F_{c''}\},
\] where
\begin{align*}
\iota_1\colon F_{{\rm Park}(a')}&\rightarrow F_{a'}\\[5pt]
\iota_2\colon F_{{\rm Park}(b')}&\rightarrow F_{b'} \\[5pt]
\iota_3\colon F_{{\rm Park}(a'')}&\rightarrow F_{a''} \\[5pt]
\iota_4\colon F_{{\rm Park}(b'')}&\rightarrow F_{b''}
\end{align*}
are the natural one-to-one correspondences   as mentioned before. It is not difficult to verify that $(\Theta,c',c'')\mapsto (a',a'',b',b'',\Theta',\Theta'')$ is a bijection. Moreover, according to our construction, we have
\[
c'=a'\Theta_{m} b',\ \ {\rm and}\ \ c''=a''\Theta_{m}b''.
\] Then equations \eqref{fact1} and \eqref{fact2} follow immediately from Lemma \ref{Parkized lemm}.
  \qed

\begin{prop}\label{algebraic map DELTA}
For any $a\in P_m$ and $b\in P_n$, we have $\Delta(M_a\star M_b)=\Delta(M_a)\star\Delta(M_b)$.
\end{prop}
\pf  By the definitions of product and coproduct given by
\eqref{product of KP} and \eqref{coproduct of KP}, we have
\begin{align*}
\Delta(M_a\star M_b)&=\Delta(\sum_{\Theta\in R(a,b)}M_{a\Theta b})
=\sum_{\Theta\in R(a,b)}\Delta(M_{a\Theta b})\\[5pt]
&=\sum_{\Theta\in R(a,b)}\sum_{F_{c'}\dotcup F_{c''}=F_{a\Theta b}}M_{{\rm Park}(c')}
\otimes M_{{\rm Park}(c'')}\\[5pt]
&=\sum_{(\Theta,c',c'')\in X(a,b)} M_{{\rm Park}(c')}\otimes M_{{\rm Park}(c'')}
\end{align*}
and
\begin{align*}
\Delta(M_a)\star\Delta(M_b)&=(\sum_{F_{a'}\dotcup F_{a''}=F_a}M_{{\rm Park}(a')}
\otimes M_{{\rm Park}(a'')})\star (\sum_{F_{b'}\dotcup F_{b''}=F_b}M_{{\rm Park}(b')}
\otimes M_{{\rm Park}(b'')})\\[5pt]
&=\sum_{F_{a'}\dotcup F_{a''}=F_a}\sum_{F_{b'}\dotcup F_{b''}=F_b} M_{{\rm Park}(a')}\star M_{{\rm Park}(b')}\otimes M_{{\rm Park}(a'')}\star M_{{\rm Park}(b'')}\\[5pt]
&=\sum_{(a',a'',b',b'',\Theta',\Theta'')\in Y(a,b)}M_{{\rm Park}(a')\Theta'{\rm Park}(b')}
\otimes M_{{\rm Park}(a'')\Theta''{\rm Park}(b'')}.
\end{align*}
By Lemma \ref{main lemm}, we deduce that $\Delta(M_a\star M_b)=\Delta(M_a)\star
\Delta(M_b)$. \qed

Now Theorem \ref{Main theorem} follows from Proposition \ref{algebra},
Proposition \ref{coalgebra} and Proposition \ref{algebraic map DELTA}.

\section{Freeness of PFSym and the $Q$-basis}\label{sec:4}
The algebra NCSym was originally studied by Wolf in \cite{Wolf1936}, where she aimed
to show that this algebra is freely generated. A combinatorial
description of the generating set of Wolf has been found by Bergeron, Reutenauer,
Rosas, and Zabrocki \cite{Bergeron-R-R-Z2008} in terms of monomial basis. The freeness of NCSym can also be proved by using another two bases, the power sum basis introduced by Rosas and Sagan \cite{Rosas-Sagan2006} and the q-basis introduced by Bergeron and Zabrocki \cite{Bergeron-Zabrocki2005}.

In this section, we shall show the freeness of the algebra  ${\rm PFSym}$. To this end, we introduce the notions of atomic parking function and unsplitable parking function.  We  prove that ${\rm PFSym}$ is freely generated by
those monomial basis elements indexed by unsplitable parking functions. Using a partial order on parking functions, we introduce a new basis $\{Q_a\}$ which is related to the monomial basis $\{M_a\}$ via M\"{o}bius inversion. In terms of this basis, we find another free generating set of PFSym indexed by atomic parking functions.
Based on this fact, we deduce that  PFSym is isomorphic to the graded dual PQSym$^{*}$ of the Hopf algebra PQSym in \cite{Hivert-Novelli-Thibon2008}.

\begin{defi}
Suppose that $a=a_1a_2\cdots a_m\in P_m$ and $b=b_1b_2\cdots b_n\in P_n$.
Let $F_a=(u_1,u_2,\ldots,u_r)$, $F_b=(w_1,w_2,\ldots,w_s)$ and let
 \begin{align}\label{split prod3}
F_{a}\circ F_{b}=\begin{cases} (u_1\cdot(w_1+m),\,\ldots,\,u_r\cdot
(w_r+m),\, w_{r+1}+m,\,\ldots,\,w_{s}+m),
&\textrm{if } r\leq s; \\[5pt]
(u_1\cdot (w_1+m),\,\ldots,\,u_s\cdot (w_s+m),\,u_{s+1},\,
\ldots,\,u_{r}), &\textrm{if } r>s.
\end{cases}
\end{align}
 By Remark \ref{remark for F_a}, we have $F_{a}\circ F_{b}\in\mathcal{F}_{m+n}$. So there is a unique parking function, denoted by $a\circ b$, such that $F_{a\circ b}=F_a\circ F_b$. We call $a\circ b$ the split product of $a$ and $b$.
\end{defi}

For example, let $a=2213\in P_4$ and $b=32214\in P_5$. Then we have $F_a=(13,22)$, $F_b+4=(58,66,7)$, and $F_a\circ F_b=(1358,2266,7)$. So we have $a\circ b=722661358.$
It  can be observed easily that the binary operation $\circ$ is associative and so we can define the split product of $k\geq 2$ parking functions:
\[
a^{(1)}\circ a^{(2)}\circ\cdots\circ a^{(k)}
\] in the obvious way.

\begin{defi}
 A parking function $a$ is said to be unsplitable if it can not be written as the split
product of two nonempty parking functions.
\end{defi}
Denote by $\mathcal{U}P_n$ the set of
unsplitable  parking functions of length $n$. For example,  we have
\[
\mathcal{U}P_3:=\{111, 112, 121, 131, 132, 211,
 212, 311,  221, 231, 321\}.
\]
Set $\mathcal{U}P=\cup_{n\geq 1}\mathcal{U}P_n$.

\begin{theo}\label{unsplitable generating set for KP}
The algebra ${\rm PFSym}$ is freely generated by the set
$\{M_a\ |\ a\in\mathcal{U}P\}$.
\end{theo}
\pf  Define a total order $\preceq^{*}$ on the set of  words on positive integers as follows. Suppose that $a=a_1a_2\cdots a_i$ and $b=b_1b_2\cdots b_j$, then $a\prec^{*} b$ if and only if
 \begin{itemize}
 \item[(1)]  $b$ is a proper  initial subword of $a$, that is, $i>j$ and $a_l=b_l$ for any $1\leq l\leq j$ or
 \item[(2)]  $a$ is not a proper initial subword of $b$ and $a<_{{\rm lex}}b$, where $<_{{\rm lex}}$ denotes the lexicographical order.
 \end{itemize}
For $a,b\in P_n$, define $a\preceq_{{\rm lex}}^{*} b$ if $F_a \preceq_{{\rm lex}}^{*} F_b$, where  $\preceq_{{\rm lex}}^{*}$ denotes the lexicographic order  induced by $\preceq^{*}$. That is $(u_1,u_2,\ldots,u_r)\prec_{{\rm lex}}^{*} (w_1,w_2,\ldots,w_s)$ if
 there exists $i\leq {\rm min}(r,s)$ such that  $u_j= w_j$ for any $1\leq j< i$ and $u_i\prec^{*} w_i$.  For example, we have
\begin{align*}
&(111)\prec_{{\rm lex}}^{*}(112)\prec_{{\rm lex}}^{*}(113)\prec_{{\rm lex}}^{*}(11,2)\prec_{{\rm lex}}^{*}(11,3)\prec_{{\rm lex}}^{*}(121)\prec_{{\rm lex}}^{*}(122)\prec_{{\rm lex}}^{*}(123)\prec_{{\rm lex}}^{*}\\[5pt]
&(12,2)\prec_{{\rm lex}}^{*}(12,3)\prec_{{\rm lex}}^{*}(131)\prec_{{\rm lex}}^{*}(132)\prec_{{\rm lex}}^{*}(13,2)\prec_{{\rm lex}}^{*}(1,22)\prec_{{\rm lex}}^{*}(1,23)\prec_{{\rm lex}}^{*}(1,2,3),
\end{align*}
 and so in $P_3$, we have
\begin{align*}
&111\prec_{{\rm lex}}^{*}112\prec_{{\rm lex}}^{*}113\prec_{{\rm lex}}^{*}211\prec_{{\rm
lex}}^{*}311\prec_{{\rm lex}}^{*}121\prec_{{\rm lex}}^{*}122\prec_{{\rm lex}}^{*}123\prec_{{\rm lex}}^{*}\\[5pt]
&212\prec_{{\rm lex}}^{*}312\prec_{{\rm lex}}^{*}131\prec_{{\rm lex}}^{*}132\prec_{{\rm lex}}^{*}213\prec_{{\rm lex}}^{*}221\prec_{{\rm lex}}^{*}231\prec_{{\rm lex}}^{*}321.
\end{align*}

For any parking function $a$, we can uniquely decompose $a$ into split product of unsplitable parking functions:
  $$a=a^{(1)}\circ a^{(2)}\circ\cdots\circ a^{(k)}.$$
Now define $R_a=M_{a^{(1)}}\star M_{a^{(2)}}\star\cdots\star M_{a^{(k)}}$.  We claim that
\begin{align}\label{triangularity}
R_a=M_a+\sum_{a\prec_{{\rm lex}}^{*} b}  c_bM_b.
\end{align}
The claim can be proved by induction on $k$. The case when $k=1$ is trivially true since $R_a=M_a$. Assume that $k>1$ and the claim is true for any word $a'$ that can be decomposed into split product of $k-1$ unsplitable parking functions. Suppose that
$a=a^{(1)}\circ a^{(2)}\circ\cdots\circ a^{(k)}\in P_n$ and $a'=a^{(1)}\circ a^{(2)}\circ\cdots\circ a^{(k-1)}\in P_m$. Then  $a=a'\circ a^{(k)}$ and
\begin{align*}
R_a&=R_{a'}\star M_{a^{(k)}}\\[5pt]
&=(M_{a'}+\sum_{a'\prec_{{\rm lex}}^{*} b'}c_{b'}M_{b'})\star M_{a^{(k)}}\\[5pt]
&=\sum_{\Theta\in R(a',a^{(k)})}M_{a'\Theta a^{(k)}}+\sum_{a'\prec_{{\rm lex}}^{*}b'}\ \sum_{\Theta\in R(b', a^{(k)})}c_{b'}M_{b'\Theta a^{(k)}}.
\end{align*}
To complete the proof, it suffices to show the following facts:

(1) $a'\circ a^{(k)}=a'\Lambda a^{(k)}$ for a unique $\Lambda\in R(a',a^{(k)})$;

(2) For any $\Theta\in R(a',a^{(k)})$, if $\Theta\neq\Lambda$, then $a'\circ a^{(k)}\prec_{{\rm lex}}^{*}a'\Theta a^{(k)}$;

(3) If $a'\prec_{{\rm lex}}^{*} b'$ in $P_m$, then $a'\circ a^{(k)}\prec_{{\rm lex}}^{*}b'\circ a^{(k)}$ in $P_n$.

Assume that $F_{a'}=(u_1,u_2,\ldots,u_r), F_{a^{(k)}}=(v_1,v_2,\ldots,v_s)$. Obviously, the matching $\Lambda$ with edge set
\[
\{(u_i,v_i+m)\ |\ 1\leq i\leq \min(r,s)\}
\]
is the desired matching in fact (1). For any $\Theta\in R(a',a^{(k)})$, let $F_{a'\Theta a^{(k)}}=(w_1,w_2,\ldots,w_t)$. A simple observation shows that if $i\leq r$, then $w_i=u_i$ or $w_i=u_i\cdot(v_j+m)$ for some $1\leq j\leq s$. If $\Theta\neq\Lambda,$ then there exists some $h\leq \min(r,s)$ such that
$w_h\neq u_h\cdot(v_h+m)$. We choose $h$ to be the smallest integer satisfying this property. Then $w_h=u_h$ or $w_h=u_h\cdot(v_j+m)$ for some $j>h$. In either case, we have
$u_h\cdot(v_h+m)\prec^{*} w_h$. By the definition of $\prec_{{\rm lex}}^{*}$, we deduce that $a'\circ a^{(k)}\prec_{{\rm lex}}^{*} a'\Theta a^{(k)}$. This completes the proof of fact (2). We continue to prove fact (3). If $a'\prec_{{\rm lex}}^{*} b'$ in $P_m$, then we can assume that $$F_{b'}=(u_1,u_2,\ldots,u_{h-1},u_{h}^{'},\ldots,u_{l}^{'}),$$
where $h\leq r$ and $u_h\prec^{*} u_{h}^{'}$. This implies that if $F_{a'\circ a^{(k)}}=(w_1,w_2,\ldots,w_{h-1},w_{h},\ldots,w_t)$, then $F_{b'\circ a^{(k)}}$ takes the form
\[
(w_1,w_2,\ldots,w_{h-1},w_{h}',\ldots,w_{j}'),
\]
where
\begin{align*}
w_{h}'=\left\{
         \begin{array}{ll}
           u_{h}', & \hbox{if $w_h=u_h$;} \\[5pt]
           u_{h}'\cdot(v_h+m), & \hbox{if $w_h=u_h\cdot(v_h+m)$.}
         \end{array}
       \right.
\end{align*}
In either case, it is easy to show that $w_h\prec^{*}w_{h}'$. Hence we have $a'\circ a^{(k)}\prec_{{\rm lex}}^{*} b'\circ a^{(k)}$. This proves fact (3).

Therefore
the claim \eqref{triangularity} is true. By triangularity, we deduce that the set $\cup_{n\geq 0}\{R_a\ |\ a\in P_n\}$ is also a
basis for $\mathbb{K}\mathcal{P}$. So by the definition of $R_a$, we conclude that
$$\{M_a\ |\ a {\rm\ is\ unsplitable}\}$$ freely generates the algebra ${\rm PFSym}.$
\qed

\vspace{0.5cm}

From the proof of Theorem \ref{unsplitable generating set for KP}, we
deduce that $\cup_{n\geq 0}\{R_a\ |\ a\in P_n\}$ is a
multiplicative basis of ${\rm PFSym}$ since by the definition of $R$-basis, we have $R_a\star R_b=R_{a\circ b}$.  We now introduce another multiplicative basis.   To begin with, we introduce the slash product of parking functions and the notion of atomic parking function.

\begin{defi}\label{definition of slash product}
Suppose that $a=a_1a_2\cdots a_m\in P_m$ and  $b=b_1b_2\cdots b_n\in P_n$.
Let $a|b=(b_1+m)(b_2+m)\cdots (b_n+m)a_1a_2\cdots a_m$. We call $a|b$ the slash product of $a$ and $b$.
\end{defi}

Let $F_a=(u_1,u_2,\ldots,u_r)$,
$F_b=(w_1,w_2,\ldots,w_s)$ and let
\[
F_{a}|F_{b}=(u_1,u_2,\ldots,u_r,w_1+m,w_2+m,\ldots,w_s+m).
\]
Clearly, we have $F_{a|b}=F_{a}|F_{b}$. So the slash product of parking functions can be seen as a natural generalization of the slash product of set partitions \cite{Chen-Li-Wang2011}. Note that the slash product  $|$ is also associative.

\begin{defi}
A parking function $a$ is said to be atomic if there are no nonempty
parking functions  $b$ and $c$ such that $a=b | c$.
\end{defi}

Note that the notion of atomic parking function is related to the notion of prime parking function introduced by Gessel in 1977, see\cite{Novelli-Thibon2007} for details. Let $\mathcal{A}P_n$ denote the set of atomic  parking functions of length $n$. For example, we have
$$
\mathcal{A}P_3=\{111, 211, 121, 131, 112, 212, 122, 132,
 113, 213, 123\}.
$$  Set $\mathcal{A}P=\cup_{n\geq 1}\mathcal{A}P_n$

Now we define a partial order $\leq_{*}$ on $P_n$ and then introduce a new basis $\{Q_a\}$ via M\"{o}bius inversion. Suppose that $a,b\in P_n$  and $F_a=(w_1,w_2,\ldots,w_r)$. We say $b$ covers $a$ if
 there exist $1\leq i<j\leq r$ such that every element in $w_i$ is less than or
equal to every element in $w_j$ and
$$F_b=(w_1,\ldots,w_{i-1},w_{i}\cdot w_j,w_{i+1},\ldots,w_{j-1},w_{j+1},\ldots,w_{r}).$$
Let $\leq_{*}$ denote the partial order on $P_n$ generated by these covering relations.
For example, the following figure shows the Hasse diagram of $(P_3,\leq_{*})$.

\begin{center}\setlength{\unitlength}{1mm}
\begin{picture}(80,50)
\put(-15,30){\circle*{1}}\put(-18,35){$111$}
\put(-3,30){\circle*{1}}\put(-6,35){$121$}
\put(9,30){\circle*{1}}\put(6,35){$131$}
\put(21,35){\circle*{1}}\put(19,37){$112$}
\put(21,25){\circle*{1}}\put(19,21){$211$}
\put(21,25){\line(0,1){10}}
\put(33,35){\circle*{1}}\put(31,37){$113$}
\put(33,25){\circle*{1}}\put(31,21){$311$}
\put(33,25){\line(0,1){10}}
\put(46,35){\circle*{1}}\put(44,37){$122$}
\put(42,25){\circle*{1}}\put(40,21){$221$}
\put(50,25){\circle*{1}}\put(48,21){$212$}
\put(42,25){\line(2,5){4}}\put(50,25){\line(-2,5){4}}
\put(68,40){\circle*{1}}\put(66,42){$123$}
\put(63,30){\circle*{1}}\put(73,30){\circle*{1}}\put(86,30){\circle*{1}}
\put(56,28){$231$}\put(74,28){$312$}\put(88,28){$213$}\put(102,28){$132$}
\put(63,30){\line(1,2){5}}\put(73,30){\line(-1,2){5}}\put(100,30){\circle*{1}}
\put(73,20){\circle*{1}}
\put(74,18){$321$}
\put(73,20){\line(-1,1){10}}\put(73,20){\line(0,1){10}}\put(73,20){\line(4,3){13}}
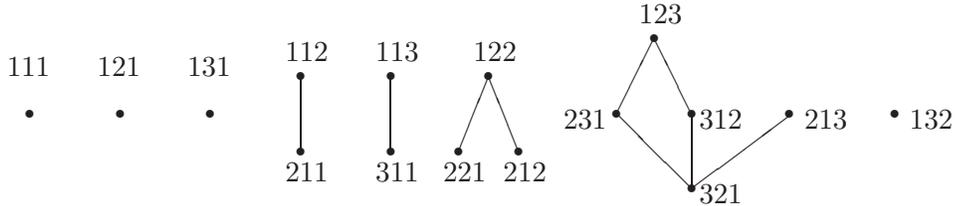
\put(18,5){Fig.1. the poset  $(P_3,\leq_{*})$}
\end{picture}
\end{center}

For any parking function $a\in P_n$, set
\begin{align}\label{new basis}
Q_a=\sum_{a\leq_{*} b}M_b.
\end{align}
By M\"{o}bius inversion, $\cup_{n\geq 0}\{Q_a\ |\ a \in P_n\}$
is also a  basis of $\mathbb{K}\mathcal{P}$. Moreover, the following theorem shows that it is a multiplicative basis.

\begin{theo} \label{atomic generating set for KP}
Suppose that $a\in P_m$ and $b\in P_n$. Then
\begin{align}\label{product for Q basis}
Q_a\star Q_b=Q_{a|b},
\end{align}
and
\begin{align}\label{coproduct for Q basis}
\Delta(Q_a)=\sum_{F_{a'}\dotcup F_{a''}=F_a}Q_{{\rm Park}(a')}\otimes Q_{{\rm Park}(a'')}.
\end{align}In particular, the set $\{Q_a\ |\ a\in\mathcal{A}P\}$ freely generates the algebra ${\rm PFSym}$.
\end{theo}

\pf By the definitions of $Q_a$ and $\star$, we have
\begin{align}\label{product-Q-1}
\nonumber Q_a\star Q_b&=(\sum_{a\leq_{*} \tilde{a}} M_{\tilde{a}})\star (\sum_{b\leq_{*} \tilde{b}} M_{\tilde{b}})\\[5pt]
\nonumber &=\sum_{a\leq_{*} \tilde{a},\ b\leq_{*} \tilde{b}} M_{\tilde{a}}\star M_{\tilde{b}}\\[5pt]
&=\sum_{a\leq_* \tilde{a},\ b\leq_* \tilde{b}}\ \sum_{\Theta\in R(\tilde{a},\tilde{b})} M_{\tilde{a}\Theta \tilde{b}},
\end{align}
and
\begin{align}\label{product-Q-2}
Q_{a| b}=\sum_{a|b\leq_{*} c} M_c.
\end{align}
According to \eqref{product-Q-1} and \eqref{product-Q-2}, for any $c'\in P_{m+n}$, the coefficient of $M_{c'}$ in $Q_{a}\star Q_b$ or in $Q_{a| b}$ is $0$ or $1$. So it suffices to show that if
\[
A(a,b)=\{\tilde{a}\Theta \tilde{b}:\  \tilde{a}\in P_m,
\tilde{b}\in P_n,  a\leq_{*}\tilde{a}, b\leq_{*}\tilde{b}, \Theta\in R(\tilde{a},\tilde{b}) \}
\]
and
\[ B(a,b)=\{c:\ c\in P_{m+n},  a|b\leq_{*} c\},
\]
then $A(a,b)=B(a,b)$. Let $F_a=(u_1,u_2,\ldots,u_r)$ and $F_{b}=(w_1,w_2,\ldots,w_s)$. Suppose that $c=\tilde{a}\Theta\tilde{b}\in A(a,b)$. Since $a\leq_{*}\tilde{a}$ and $b\leq_{*} \tilde{b}$, we deduce that each word
$v$ in $F_{\tilde{a}}$ takes the form
\[
v=u_{i_1}\cdot u_{i_2}\cdots u_{i_k},\   1\leq i_1<i_2<\cdots <i_k\leq r
\]
and each word $v'$ in $F_{\tilde{b}}$ takes the form
\[
v'=w_{j_1}\cdot w_{j_2}\cdots w_{j_l},\   1\leq j_1<j_2<\cdots <j_l\leq s
\]
where each integer in $u_{i_{t}}$ is less than or equal to each integer in $u_{i_{t+1}}$, and each integer in $w_{j_{t}}$ is less than or equal to each integer in $w_{j_{t+1}}$.
So each word $v^{*}$ in $F_c$ takes one of the following three forms:\\[5pt]
$(1)\ \ \ v^{*}=u_{i_1}\cdot u_{i_2}\cdots u_{i_k}$\\[5pt]
$(2)\ \ \  v^{*}=(w_{j_1}+m)\cdot (w_{j_2}+m)\cdots (w_{j_l}+m)$\\[5pt]
$(3)\ \ \  v^{*}=u_{i_1}\cdot u_{i_2}\cdots u_{i_k}\cdot (w_{j_1}+m)\cdot (w_{j_2}+m)\cdots (w_{j_l}+m).$\\[5pt]
By the definition of $\leq_{*}$, we have $a|b\leq_{*} c$ and so $c\in B(a,b)$. Conversely, if $c\in B(a,b)$, then each word $v^{*}$ in $F_{c}$ takes one of the above three forms.
Let $\tilde{a},\tilde{b}$ be the parking functions determined by
\begin{align*}
F_{\tilde{a}}&=\{v\ |\ v\in F_c {\rm\ and}\ \max(v)\leq m\}\cup \{v\ |\ v\cdot (v'+m)\in F_c {\rm\ for\ some}\ v'\ {\rm\ and}\ \max(v)\leq m\}\\[5pt]
F_{\tilde{b}}&=\{v'\ |\ v'+m\in F_c \}\cup \{v'\ |\ v\cdot(v'+m)\in F_c {\rm\ for\ some}\ v\ {\rm\ such\ that}\ \max(v)\leq m\}.
\end{align*}
Clearly, we have $a\leq_{*}\tilde{a}$, $b\leq_{*}\tilde{b}$ and $c=\tilde{a}\Theta\tilde{b}$, where $\Theta$ is the matching between $F_{\tilde{a}}$ and $F_{\tilde{b}}+m$ with edge set
\[
\{(v,v'+m)\ |\ v\cdot (v'+m)\in F_c\}.
\]So we have $c\in A(a,b)$. Hence, $A(a,b)=B(a,b)$ and the equation \eqref{product for Q basis} holds.

By the definition of $Q_a$, we have
\begin{align*}
\Delta(Q_a)&=\Delta(\sum_{a\leq_{*} b}M_b)=\sum_{a\leq_{*}b}\Delta(M_b)\\[5pt]
&=\sum_{a\leq_{*}b}\ \sum_{F_{b'}\dotcup F_{b''}=F_b} M_{{\rm Park}(b')}\otimes M_{{\rm Park}(b'')}
\end{align*}
and
\begin{align*}
\sum_{F_{a'}\dotcup {F_{a''}=F_a}} Q_{{\rm Park}(a')}\otimes Q_{{\rm Park}(a'')}=\sum_{F_{a'}\dotcup F_{a''}=F_a}\ \sum_{{\rm Park}(a')\leq_*c'\atop
{\rm Park}(a'')\leq_* c''}M_{c'}\otimes M_{c''}.
\end{align*}
Let
\[
X(a):=\{(b,b',b'')\ |\ a\leq_* b, F_{b'}\dotcup F_{b''}=F_b\}
\]
and
\[
Y(a):=\{(a',a'',c',c'')\ |\ F_{a'}\dotcup F_{a''}=F_a, {\rm Park}(a')\leq_* c', {\rm Park}(a'')\leq_* c''\}.
\]
To show \eqref{coproduct for Q basis}, it suffices to find a bijection $\phi: (b,b',b'')\rightarrow (a',a'',c',c'')$ between $X(a)$ and $Y(a)$ such that
\begin{align*}
{\rm Park}(b')=c',\ \ {\rm Park}(b'')=c''.
\end{align*}
Let $F_{a}=(u_1,u_2,\ldots,u_r)$ and $(b,b',b'')\in X(a)$. Since $a\leq_* b$, we deduce that every word $v\in F_b$ has the form
\[
v=u_{i_1}\cdot u_{i_2}\cdots u_{i_k}, \ \ 1\leq i_1<i_2<\cdots<i_k\leq r
\] where each integer in $u_{i_t}$ is less than or equal to each integer in $u_{i_{t+1}}$. Let
\[
U(v)=\{u_{i_1},u_{i_2},\ldots,u_{i_k}\}
\]
and  let $a', a''$ be the subwords of $a$ such that
$$
F_{a'}=\bigcup_{v\in F_{b'}}U(v),\ \ F_{a''}=\bigcup_{v\in F_{b''}}U(v).
$$Clearly, $F_{a'}\cap F_{a''}=\emptyset$ and $F_{a'}\cup F_{a''}=F_a$.
Recall that there are natural bijections $\iota_1: F_{\rm Park(a')}\rightarrow F_{a'}$ and $\iota_2: F_{{\rm Park}(a'')}\rightarrow F_{a''}$.
For $v=u_{i_1}\cdot u_{i_2}\cdots u_{i_k}\in F_{b}$, we let
\begin{align*}
\iota(v)=\left\{
             \begin{array}{ll}
               \iota_1^{-1}(u_{i_1})\cdot \iota_1^{-1}(u_{i_2})\cdots \iota_1^{-1}(u_{i_k}), & \hbox{if $v\in F_{b'}$;} \\[7pt]
               \iota_2^{-1}(u_{i_1})\cdot \iota_2^{-1}(u_{i_2})\cdots \iota_2^{-1}(u_{i_k}), & \hbox{if $v\in F_{b''}$,}
             \end{array}
           \right.
\end{align*}
and let $c',c''$ be the parking functions of length $l(a')$ and $l(a'')$ respectively such that
\[
F_{c'}=\{\iota(v)\ | v\in F_{b'}\},\ \ F_{c''}=\{\iota(v)\ |\ v\in F_{b''}\}.
\] Keeping in mind that the operator ${\rm Park}$ doesn't change relative order of integers in a word, we can easily verify that $(a',a'',c',c'')\in Y(a)$ and $\phi:(b,b',b'')\rightarrow (a',a'',c',c'')$ is the desired bijection. This completes the proof of equation \eqref{coproduct for Q basis}.

Since each parking function $a\in P$ can be uniquely decomposed into slash product of atomic parking functions, it follows from \eqref{product for Q basis} that the algebra ${\rm PFSym}$ is freely generated by $\{Q_a\ |\ a\in\mathcal{A}P\}$.
\qed

%
%

Let $U(L(X))$ denote the universal enveloping algebra of  the free Lie algebra $L(X)$ on the set $X$. By Theorem \ref{unsplitable generating set for KP}, Theorem \ref{atomic generating set for KP} and  the Milnor-Moore theorem (see \cite[p. 244]{Milnor-Moore1965} or \cite[p. 274]{Sweedler1969}),  we find that
\[
{\rm PFSym}\cong U(L(\mathcal{A}P))\cong U(L(\mathcal{U}P)).
\]
To establish a more explicit isomorphism, one can use the method given by Lauve and Mastnak \cite{Lauve-Mastnak2011} to find the algebraically independent generators of the
Lie algebra of primitive elements of ${\rm PFSym}$.

Recall that a parking function $a\in P_n$ is connected if there are no nonempty parking functions $b$ and $c$ such that
$$a= b\cdot (c+l(b)),
$$
  where $l(b)$ denotes the length of $b$. Let $\mathcal{C}P$ be the set of all connected parking functions and let PQSym$^{*}$ be the graded dual of the commutative Hopf algebra PQSym introduced by Hivert, Novelli and Thibon \cite{Hivert-Novelli-Thibon2008}.  By a similar argument as \cite[Theorem 2.5]{Hivert-Novelli-Thibon2008}, we have  ${\rm PQSym}^{*}\cong U(L(\mathcal{C}P))$. Therefore, by the obvious fact that
 $a=a_1a_2\cdots a_n\in\mathcal{C}P$  if and only if $\hat{a}=a_na_{n-1}\cdots a_1\in\mathcal{A}P$, we conclude that ${\rm PFSym}\cong {\rm PQSym}^{*}$.

To conclude this section, we remark that the Hopf algebra ${\rm PFSym}$ is also cofree. In fact, since each parking function $a$ can be uniquely decomposed into slash product of atomic parking functions:
\[
a=a^{(1)}|a^{(2)}|\cdots |a^{(r)},
\]
 we can identify  parking  functions with  words on
$\mathcal{A}P$. By rephrasing the proof given in \cite[Section 4]{Bergeron-Zabrocki2005},
 one can show that ${\rm PFSym}^{*}$ is isomorphic to a shuffle algebra on a vector space with a basis indexed by atomic parking functions. Therefore the freeness of ${\rm PFSym}^{*}$, and so the cofreeness of ${\rm PFSym}$, follows from \cite[Theorem 6.1]{Reutenauer1993}.

\section{Hopf subalgebras of ${\rm PFSym}$}\label{Hopf subalgebras}\label{sec:5}
In this section, we discuss Hopf subalgebras of PFSym.  As will be seen, the Hopf algebra
${\rm NCSym}$ can be embedded as a Hopf subalgebra into PFSym in a natural way. Moreover, when restricting to permutations and non-increasing parking functions, we find two Hopf subalgebras isomorphic to the Grossman-Larson Hopf algebras of ordered trees and heap-ordered trees respectively.

For each $n\geq 0$, let
\begin{align*}
N_n=\{a\in P_n &|\ {\rm each\ word\ in}\ F_a\ {\rm is\ nondecreasing}\ \}\\[5pt]
D_n=\{a\in P_n &|\ F_a=(w_1,w_2,\ldots,w_k)\ {\rm and}\ w_i\cap w_j=\emptyset\}\\[5pt]
\mathfrak{S}_n=\{a\in P_n &|\ a\ {\rm is\ a\ permutation\ of\ }[n]\}\\[5pt]
C_n=\{a\in P_n &|\ a=a_1a_2\cdots a_n\ {\rm with}\ a_1\geq a_2\geq\cdots\geq a_n\}
\end{align*}
where we use $w_i\cap w_j$ to denote the
 intersection of the underlying sets of $w_i$  and $w_j$. Then we have
\[
\mathfrak{S}_n\subseteq D_n\
 \  {\rm and}\ \  C_n\subseteq N_n\cap D_n.
\]
Moreover, it is easy to check that each $a\in C_n$ is a minimal element in the poset
$(P_n,\leq_{*})$ and
\begin{align}\label{minimal elements}
\{b\in P_n |\ b\geq_{*} a\ {\rm for\ some}
\ a\in C_n\}=N_n\cap D_n.
\end{align}

Now set
$$
\begin{array}{cccc}
 & \mathcal{A}C_n=\mathcal{A}P_n\cap C_n,
 & \mathcal{A}C=\cup_{n\geq 1}\mathcal{A}C_n,& \  \\[5pt]
   \mathcal{A}N_n=\mathcal{A}P_n\cap N_n, & \mathcal{U}N_n=\mathcal{U}P_n\cap N_n,
 & \mathcal{A}N=\cup_{n\geq 1}\mathcal{A}N_n,& \mathcal{U}N=\cup_{n\geq 1}\mathcal{U}N_n,\\[5pt]
   \mathcal{A}D_n=\mathcal{A}P_n\cap D_n, & \mathcal{U}D_n=\mathcal{U}P_n\cap D_n,
 & \mathcal{A}D=\cup_{n\geq 1}\mathcal{A}D_n,& \mathcal{U}D=\cup_{n\geq 1}\mathcal{U}D_n,\\[5pt]
 \mathcal{A}\mathfrak{S}_n=\mathcal{A}P_n\cap \mathfrak{S}_n, & \mathcal{U}\mathfrak{S}_n=\mathcal{U}P_n\cap \mathfrak{S}_n,
 & \mathcal{A}\mathfrak{S}=\cup_{n\geq 1}\mathcal{A}\mathfrak{S}_n,& \mathcal{U}\mathfrak{S}=\cup_{n\geq 1}\mathcal{U}\mathfrak{S}_n,

\end{array}
$$
and set
\begin{align*}
\mathbb{K}\mathcal{N}&=\bigoplus_{n\geq 0}\mathbb{K}\{M_a\ |\ a\in N_n\},\qquad
\mathbb{K}\mathcal{D}=\bigoplus_{n\geq 0}\mathbb{K}\{M_a\ |\ a\in D_n\},\\[5pt]
\mathbb{K}\mathfrak{S}&=\bigoplus_{n\geq 0}\mathbb{K}\{M_a\ |\ a\in \mathfrak{S}_n\},\qquad
\mathbb{K}\mathcal{C}=\bigoplus_{n\geq 0}\mathbb{K}\{Q_a\ |\ a\in C_n\}.
\end{align*}

\begin{theo}\label{Hopf subalgebras}
All the subspaces $\mathbb{K}\mathcal{N}$, $\mathbb{K}\mathcal{D}$,
$\mathbb{K}\mathfrak{S}$ and $\mathbb{K}\mathcal{C}$ are Hopf
subalgebras of {\rm PFSym}. Moreover, we have
\begin{itemize}
\item[(1)] the Hopf algebra $\mathbb{K}\mathcal{N}$ has two free generating sets:
$\{M_a\ |\ a\in\mathcal{U}N\}$ and $\{Q_a\ |\ a\in\mathcal{A}N\}$;
\item[(2)] the Hopf algebra $\mathbb{K}\mathcal{D}$ has two free generating sets:
$\{M_a\ |\ a\in\mathcal{U}D\}$ and $\{Q_a\ |\ a\in\mathcal{A}D\}$;
\item[(3)] the Hopf algebra $\mathbb{K}\mathfrak{S}$ has two free generating sets:
$\{M_a\ |\ a\in\mathcal{U}\mathfrak{S}\}$ and $\{Q_a\ |\ a\in\mathcal{A}\mathfrak{S}\}$;
\item[(4)] the Hopf algebra $\mathbb{K}\mathcal{C}$ is freely  generated by
$\{Q_a\ |\ a\in\mathcal{A}C\}$.
\end{itemize}
\end{theo}
\pf Suppose that $a\in N_m$ (resp. $D_m$, $\mathfrak{S}_m$) and $b\in N_n$
(resp. $D_n$, $\mathfrak{S}_n$). It is easy to show that
\[
a\Theta b\in N_{m+n} \ \ ({\rm resp.}\ D_{m+n}, \mathfrak{S}_{m+n})
\]
for any $\Theta\in R(a,b)$. So the product is closed on $\mathbb{K}\mathcal{N}$ (resp.
$\mathbb{K}\mathcal{D}$, $\mathbb{K}\mathfrak{S}$) . Since the operator
${\rm Park}$ preserve the relative order of integers in words, the coproduct is also closed on
$\mathbb{K}\mathcal{N}$ (resp. $\mathbb{K}\mathcal{D}$, $\mathbb{K}\mathfrak{S}$). Now we consider
 the space $\mathbb{K}\mathcal{C}$. Note that $C_n$ consists of parking functions $a\in P_n$ such
that each word in $F_a$ is a word on a singleton set. So for any
 $a\in C_m, b\in C_n$ and subword $a'$ of $a$ such that $F_{a'}\subseteq F_a$, we have
\[
a|b\in C_{m+n}\ \ {\rm and}\ \ \
{\rm Park}(a')\in C_k\ \ \ \
 {\rm for\  some}\  k\leq m.\]
Hence it follows from
 \eqref{product for Q basis} and \eqref{coproduct for Q basis} that
\[
\mathbb{K}\mathcal{C}\cdot \mathbb{K}\mathcal{C}\subseteq \mathbb{K}\mathcal{C}
\] and
$$\Delta(\mathbb{K}\mathcal{C})\subseteq \mathbb{K}\mathcal{C}\otimes \mathbb{K}\mathcal{C}.$$
So   $\mathbb{K}\mathcal{N}$, $\mathbb{K}\mathcal{D}$,
$\mathbb{K}\mathfrak{S}$ and $\mathbb{K}\mathcal{C}$ are Hopf
subalgebras of PFSym.

Note that for each $n\geq 1$, the set $N_n$ is a dual order ideal of $(P_n,\leq_{*})$, i.e., if $a\in N_n$ ,$b\in P_n$ and $b\geq_{*} a$, then  $b\in N_n$. So we have
$\cup_{n\geq 0}\{Q_a\ |\ a\in N_n\}$  is also a basis for $\mathbb{K}\mathcal{N}$.
Now the statement (1) follows from the proof  of
Theorem \ref{unsplitable generating set for KP} and equation \eqref{product for Q basis}. Similarly, we can prove the statements (2) and (3).
Using equation \eqref{product for Q basis} again, we get
the proof of statement (4).
\qed

It should be noted  that $C_n$ is not a dual order ideal of $(P_n,\leq_{*})$ in general. In fact, by \eqref{minimal elements}, $C_n$ consists of minimal elements of the dual
order  ideal $N_n\cap D_n$.
So $M_a\not\in\mathbb{K}\mathcal{C}$ for any $a\in C_n$  unless $a=11\cdots 1$.

\begin{coro}\label{embedding omega}
Let $\tilde{\Pi}_n=N_n\cap \mathfrak{S}_n$. Then
$$\mathbb{K}\tilde{\Pi}=\bigoplus_{n\geq 0} \mathbb{K}\{M_a\ |\ a\in\tilde{\Pi}_n\}$$
is a Hopf subalgebra of ${\rm PFSym}$ which is isomorphic
to the Hopf algebra ${\rm NCSym}$.
\end{coro}

\pf By Theorem \ref{Hopf subalgebras}, $\mathbb{K}\tilde{\Pi}$ is a Hopf subalgebra of
 ${\rm PFSym}$. We proceed to construct an explicit isomorphism between NCSym and
$\mathbb{K}\tilde{\Pi}$.
For any finite set of integers $B$, let $w(B)$ denote the increasing
arrangement of elements in $B$. Suppose that $\pi=\{B_1,B_2,\cdots,B_k\}\in\Pi_n$. Let $\omega(\pi)=w(B_k)\cdot w(B_{k-1})\cdots w(B_1)$,
where $\cdot$ denotes concatenation of words. Set $\omega(\emptyset)=\epsilon$. It can
be checked easily that,  for each $n\geq 0$, $\omega: \Pi_n\rightarrow \tilde{\Pi}_n$ is a bijection.
Now define a linear map
$\bar{\omega}: {\rm NCSym}\rightarrow \mathbb{K}\tilde{\Pi}$ by setting
$\bar{\omega}(M_\pi)=M_{\omega(\pi)}$. It is straightforward to show that $\bar{\omega}$ is a Hopf algebra isomorphism.
\qed

Note that the subposet $(\tilde{\Pi}_n,\leq_{*})$ of $(P_n,\leq_{*})$ is isomorphic to the poset
$(\Pi_n,\leq_{*})$  introduced by Bergeron and Zabrocki \cite{Bergeron-Zabrocki2005}. As a matter of fact,  it is easy to see  that the bijection
$\omega:\Pi_n\rightarrow \tilde{\Pi}_n$  given in the proof of Corollary \ref{embedding omega} is order-preserving. Therefore,
the  basis $\{Q_a\ |\ a\in \cup_{n\geq 0}\tilde{\Pi}_n\}$  for $\mathbb{K}\tilde{\Pi}$
 corresponds to the basis
$\{{\bf q}_{\pi}\ |\ \pi\in \cup_{n\geq 0}\Pi_n\}$
 for ${\rm NCSym}$.

 By Theorem \ref{Hopf subalgebras} and the Milnor-Moore theorem, we find that
\begin{align*}
&\mathbb{K}\mathcal{N}\cong U(L(\mathcal{A}\mathcal{N}))\ \qquad
\mathbb{K}\mathcal{D}\cong U(L(\mathcal{A}\mathcal{D})),\\[5pt]
&\mathbb{K}\mathfrak{S}\cong U(L(\mathcal{A}\mathfrak{S}))\ \qquad\
\mathbb{K}\mathcal{C}\cong U(L(\mathcal{A}\mathcal{C})),
\end{align*}
where $L(X)$ denotes the free Lie algebra on the set $X$ and  $U(L)$ denotes the universal
enveloping algebra of  $L$. The Hopf algebras $\mathbb{K}\mathcal{N}$ and $\mathbb{K}\mathcal{D}$ seem to be new to our knowledge, whereas the Hopf algebras  $\mathbb{K}\mathfrak{S}$ and $\mathbb{K}\mathcal{C}$  are well-known  combinatorial Hopf algebras
studied in different ways from the literatures.

\begin{coro}
The Hopf algebra $\mathbb{K}\mathfrak{S}$ is isomorphic to the Grossman-Larson Hopf algebra $\mathcal{H}_{HO}$ of heap-ordered trees.
\end{coro}
\pf By Theorem \ref{Hopf subalgebras} and the Milnor-Moore theorem, we find that
\[
\mathbb{K}\mathfrak{S}\cong U(L(\mathcal{A}\mathfrak{S})).
\]
Let $\mathfrak{S}{\bf QSym}^{*}$ denote the graded dual of the Hopf algebra $\mathfrak{S}{\rm QSym}$ in \cite{Hivert-Novelli-Thibon2008}. By \cite[Theorem 3.4]{Hivert-Novelli-Thibon2008}), we have
\[
\mathfrak{S}{\bf QSym}^{*}\cong U(L(X)),
\]
where
\[
X=\cup_{n\geq 1}\{ \pi\in\mathfrak{S}_n\ |\ \pi_1\pi_2\cdots\pi_i\notin\mathfrak{S}_i\ {\rm for\ any\ }1\leq i\leq n-1\},
\]denotes the set of  connected permutations. Since  a permutation $a_1a_2\cdots a_n$ is atomic if and only if $a_{n}a_{n-1}\cdots a_1$ is connected, we deduce that
\[
\mathbb{K}\mathfrak{S}\cong \mathfrak{S}{\bf QSym}^{*}.
\]
Now the conclusion follows from \cite[Corollary 3.5]{Hivert-Novelli-Thibon2008}, which states that $\mathfrak{S}{\bf QSym}^{*}\cong \mathcal{H}_{HO}$.
\qed

We remark   that there are another two Hopf algebras $\Phi{\bf Sym}$ (see \cite[Proposition 4.3]{Hivert-Novelli-Thibon2008}) and $\mathbb{K}\tilde{\mathfrak{S}}$ (see \cite{Grossman-Larson2009}, where the authors use the notation $\mathbb{K}\mathfrak{S}$ we use here) built on permutations  which are isomorphic to the Hopf algebra $\mathcal{H}_{HO}$. Their coalgebra structures are defined in a similar way to that of $\mathbb{K}\mathfrak{S}$.
More specifically, for a permutation $a$ with $F_a=\{w_1,w_2,\ldots,w_k\}$, let $\varphi(a)$ be the permutation $(w_1)(w_2)\cdots (w_k)$ which is expressed as a product of cycles. Then the map $\varphi$ induces a coalgebra isomorphism $M_a \mapsto \varphi(a)$ from $\mathbb{K}\mathfrak{S}$ to $\Phi{\bf Sym}$ or $\mathbb{K}\tilde{\mathfrak{S}}$.   However, $\varphi$ is not an algebra map since the product of $\mathbb{K}\mathfrak{S}$ is different from both the product given by cyclic shuffle associated with matchings on cycles and the heap product by Grossman and Larson.

  As for the Hopf algebra $\mathbb{K}\mathcal{C}$, we have
\begin{coro}
The Hopf algebra $\mathbb{K}\mathcal{C}$ is isomorphic to the Grossman-Larson Hopf algebra $\mathcal{H}_{O}$ of ordered trees.
\end{coro}

\pf Recall that  the number $|C_n|$ of non-increasing parking functions of
length $n$ is given by the $n$th Catalan number. Since it can be shown easily that
\[
\mathcal{A}C_n=\{u\cdot 1\ |\ u\in C_{n-1}\},
\]
 it follows that the set $\mathcal{A}C_n$ is enumerated by the $(n-1)$th Catalan number.
So we deduce that,  as a free  algebra generated by the Catalan set
$\{Q_a\ |\ a\in\mathcal{A}C\}$,   the Hopf algebra $\mathbb{K}\mathcal{C}$  must be
isomorphic to the Grossman-Larson Hopf algebra  $\mathcal{H}_{O}$ of ordered trees \cite{Grossman-Larson1989}.
\qed

The Hopf algebra $\mathcal{H}_{O}$ also appears as the graded dual of the quotient of the Loday-Ronco Hopf algebra  of planar binary trees by  its coradical filtration
\cite{Aguiar-Sottile2005-2,Aguiar-Sottile2006,Loday-Ronco1998}. Besides, this Hopf algebra is also isomorphic to the Catalan Quasi-symmetric Hopf algebra ${\rm CQSym}$, which is a   Hopf subalgebra  of the non-cocommutative Hopf algebra of parking functions introduced by Novelli and Thibon \cite{Novelli-Thibon2007}. So we have obtained a new way to approach these Hopf algebras.

\vspace{0.5cm}
 \noindent{\bf Acknowledgments.}
 This work was supported by the National Science Foundation of China (Grant No. 11326222, 11401316),  the Research Fund for the Doctoral Program of Higher Education of China (Grant No. 20130182120030), the Fundamental Research Funds for Central Universities (Grant No. XDJK2013C133), and the Southwest University of China (Grant No. SWU112040).

\end{document}